\documentclass{commat}

\title{%
    On the $p$-biharmonic submanifolds and stress $p$-bienergy tensors
    }

\author{%
    Khadidja Mouffoki and Ahmed Mohammed Cherif
    }

\affiliation{
    \address{Khadidja Mouffoki --
    University Mustapha Stambouli Mascara, Faculty of Exact Sciences, Mascara 29000, Algeria.
        }
    \email{%
    khadidja.mouffoki@univ-mascara.dz
    }
    \address{Ahmed Mohammed Cherif --
    University Mustapha Stambouli Mascara, Faculty of Exact Sciences, Mascara 29000, Algeria.
        }
    \email{%
    a.mohammedcherif@univ-mascara.dz
    }
    }

\abstract{%
    In this paper, we consider $p$-biharmonic submanifolds of a
space form. We give the necessary and
sufficient conditions for a submanifold to be $p$-biharmonic in a
space form. We present
some new properties for the stress $p$-bienergy tensor.
    }



\VOLUME{31}
\YEAR{2023}
\NUMBER{1}
\firstpage{117}
\DOI{https://doi.org/10.46298/cm.10289}

\begin{paper} % Use \paper, instead, if the compilation returns an error here.

\section{Introduction}
Consider a smooth map $\varphi:(M,g)\longrightarrow (N,h)$ between Riemannian manifolds, and let
$p\geq2$, for any compact domain $D$ of $M$ the $p$-energy functional of $\varphi$ is defined by
\begin{equation}\label{eq1.1}
E_{p}(\varphi;D)=\frac{1}{p}\int_{D}|d\varphi|^pv_{g},
\end{equation}
where $|d\varphi|$ is the Hilbert-Schmidt norm of the differential $d\varphi$, and $v^g$ is the volume element on $(M,g)$.
A map is called $p$-harmonic
if it is a critical point of the $p$-energy functional over any compact subset $D$ of $M$.
Let $\{\varphi_{t}\}_{t\in (-\epsilon,\epsilon)}$ be a smooth variation of $\varphi$ supported in $D$. Then
\begin{equation}\label{eq1.2}
    \frac{d}{dt}E_{p}(\varphi_{t};D)\Big|_{t=0}=-\int_{D}h(\tau_{p}(\varphi),v)v_{g},
\end{equation}
where $\displaystyle v=\frac{\partial \varphi_{t}}{\partial t}\Big|_{t=0}$ denotes the variation vector field of $\varphi$,
\begin{equation}\label{eq1.3}
    \tau_{p}(\varphi)
    =\operatorname{div}^M(|d\varphi|^{p-2}d\varphi).
\end{equation}
Let $\tau(\varphi)$ be the tension field of $\varphi$ defined by
\begin{equation}\label{eq1.4}
    \tau(\varphi)=\operatorname{trace}_g\nabla d\varphi
    =\sum_{i=1}^m\big\{\nabla^\varphi _{e_i}d\varphi(e_i)-d\varphi(\nabla^M _{e_i} e_i)\big\}.
\end{equation}
(see \cite{BW}), where $\{e_1, \dotsc ,e_m\}$ is an orthonormal frame on $(M,g)$, $m=\dim M$, $\nabla^{M}$ is the Levi-Civita connection of $(M,g)$,
and $\nabla^{\varphi}$ denotes the pull-back connection on $\varphi^{-1}TN$.
If $|d\varphi|_x\neq0$ for all $x\in M$, the map $\varphi$ is $p$-harmonic if and only if (see \cite{BG}, \cite{BI}, \cite{ali})
\begin{equation}\label{eq1.5}
   |d\varphi|^{p-2}\tau(\varphi)+(p-2)|d\varphi|^{p-3} d\varphi(\operatorname{grad}^M|d\varphi|)=0.
\end{equation}
\begin{example}
Let $n\geq2$. The inversion map
\begin{eqnarray*}
	\varphi:\mathbb{R}^n\backslash\{0\}  &\longrightarrow& \mathbb{R}^n\backslash\{0\},  \\
      x &\longmapsto& \frac{x}{\lvert x\rvert^l}
\end{eqnarray*}
is $p$-harmonic if and only if $l=\frac{n+p-2}{p-1}.$
\end{example}
A natural generalization of $p$-harmonic maps is given by integrating the square of the norm of $\tau_p(\varphi)$. More precisely, the $p$-bienergy functional of $\varphi$ is defined by
\begin{equation}\label{eq1.6}
    E_{2,p}(\varphi;D)=\frac{1}{2}\int_D|\tau_p(\varphi)|^2 v_g.
\end{equation}
We say that $\varphi$ is a $p$-biharmonic map if it is a critical point of the $p$-bienergy functional, that is to say, if it satisfies the Euler-Lagrange equation of the functional \eqref{eq1.6}, that is (see~\cite{cherif2})
\begin{eqnarray}\label{eq1.7}
\tau_{2,p}(\varphi)
   &=&\nonumber -|d\varphi|^{p-2}\operatorname{trace}_gR^{N}(\tau_{p}(\varphi),d\varphi)d\varphi
       -\operatorname{trace}_g\nabla^\varphi |d\varphi|^{p-2} \nabla^\varphi \tau_{p}(\varphi)\\
    &&-(p-2)\operatorname{trace}_g\nabla \langle \nabla^\varphi\tau_{p}(\varphi),d\varphi \rangle|d\varphi|^{p-4}d\varphi=0.
\end{eqnarray}
Let $\{e_1,\dotsc,e_m\}$ be an orthonormal frame on $(M,g)$, we have
$$\operatorname{trace}_gR^{N}(\tau_{p}(\varphi),d\varphi)d\varphi   = \sum_{i=1}^m R^N(\tau_{p}(\varphi),d\varphi(e_i))d\varphi(e_i), $$
$$\operatorname{trace}_g\nabla^\varphi |d\varphi|^{p-2} \nabla^\varphi \tau_{p}(\varphi)   = \sum_{i=1}^m\left(\nabla^{\varphi}_{e_i}|d\varphi|^{p-2}\nabla^{\varphi}_{e_i}\tau_{p}(\varphi)
-|d\varphi|^{p-2}\nabla^{\varphi}_{\nabla^{M}_{e_i}e_i}\tau_{p}(\varphi)\right),$$
$$\langle \nabla^\varphi\tau_{p}(\varphi),d\varphi \rangle=\sum_{i=1}^mh\left(\nabla^\varphi_{e_i}\tau_{p}(\varphi),d\varphi(e_i)\right),$$
\begin{eqnarray*}
\operatorname{trace}_g\nabla \langle \nabla^\varphi\tau_{p}(\varphi),d\varphi \rangle|d\varphi|^{p-4}d\varphi   &=& \sum_{i=1}^m\Big( \nabla^{\varphi}_{e_i} \langle \nabla^\varphi\tau_{p}(\varphi),d\varphi \rangle|d\varphi|^{p-4}d\varphi(e_i)\\
&&-\langle \nabla^\varphi\tau_{p}(\varphi),d\varphi \rangle|d\varphi|^{p-4}d\varphi(\nabla^{M}_{e_i}e_i)\Big).
\end{eqnarray*}
The $p$-energy functional (resp. $p$-bienergy functional) includes as a special case $(p = 2)$ the energy functional (resp. bi-energy functional), whose critical points are the usual harmonic maps
(resp. bi-harmonic maps), for more details on the concept of harmonic and bi-harmonic maps
see \cite{ES}, \cite{Jiang}.\\
$p$-harmonic maps are always $p$-biharmonic maps by definition. In particular, if $(M, g)$  is a compact orientable Riemannian manifold without boundary, and
$(N, h)$ is a Riemannian manifold with non-positive sectional curvature.
Then, every $p$-biharmonic map  from $(M,g)$ to $(N,h)$ is $p$-harmonic.
\begin{example}[\cite{cherif2}]
Let $M$ the manifold $\mathbb{R}^2\backslash\{(0,0)\}\times\mathbb{R}$ equipped with the Riemannian metric $g=(x_1^2+x_2^2)^{-\frac{1}{p}}(dx_1^2+dx_2^2+dx_3^2)$, and let
$N$ the manifold $\mathbb{R}^2$ equipped with the Riemannian metric $h=dy_1^2+dy_2^2$.
The map
\[
\varphi : (M,g) \to (N,h)
\qquad \textup{defined by} \qquad 
\varphi(x_1,x_2,x_3) = \left( \sqrt{x_1^2+x_2^2} ,\ x_3 \right)
\]
is proper $p$-biharmonic.
\end{example}

A submanifold in a Riemannian manifold is called a $p$-biharmonic submanifold if the isometric immersion defining the submanifold is a $p$-biharmonic map. We will call proper $p$-biharmonic submanifolds a $p$-biharmonic submanifols which is non $p$-harmonic.\\
In this paper, we will focus our attention on $p$-biharmonic submanifolds of space form, we give the necessary and sufficient conditions for submanifolds to be $p$-biharmonic. Then, we apply this general result to many particular cases. We also consider the stress $p$-bienergy tensor associated to the $p$-bienergy functional, and we give the relation between the divergence of the stress $p$-bienergy tensor and the $p$-bitension field (\ref{eq1.7}). Finally, we classify maps between Riemannian manifolds with vanishing stress $p$-bienergy tensor.
%
%
%
%
%
%
%
%
%
%
%
%
%
%%%%%%%%%%%%%%%%%%%%%%%%%%%%%%%%%%%%%%%%%%%%%%%%%%%%%%%%%%%%%%%%%%%%%%%%%%%%%%%%%%
%%%%%%%%%%%%%%%%%%%%%%%%%% Contents of Section 2 %%%%%%%%%%%%%%%%%%%%%%%%%%%%%%%%%
%%%%%%%%%%%%%%%%%%%%%%%%%%%%%%%%%%%%%%%%%%%%%%%%%%%%%%%%%%%%%%%%%%%%%%%%%%%%%%%%%%
\section{Main Results}

Let $M$ be a submanifold of space form $N(c)$ of dimension $m$, $\mathbf{i} : M \hookrightarrow N(c)$ be the canonical inclusion,
and let $\{e_1,\dotsc,e_m\}$ be an orthonormal frame with respect to induced Riemannian metric $g$ on $M$ by the inner product $\langle , \rangle$ on $N(c)$.
We denote by $\nabla^N$ (resp. $\nabla^M$) the  Levi-Civita connection of
$N^{n}(c)$ (resp. of $(M,g)$), by $\operatorname{grad}^M$ the gradient operator in $(M,g)$,
by $B$ the second fundamental form of the submanifold $(M,g)$, by $A$ the shape operator,  by $H$
the mean curvature vector field of $(M,g)$, and by $\nabla^\perp$ the normal connection of $(M,g)$ (see for example \cite{BW}).
Under the notation above we have the following results.

\begin{theorem}
 The canonical inclusion $\mathbf{i}$ is $p$-biharmonic if and only if
   \begin{equation}\label{ortho}
\left\{
\begin{array}{lll}
-\Delta^{\perp}  H + \operatorname{trace}_g  B( \cdot , A_H (\cdot)) - m\left(c-(p-2) |H|^2 \right)   H &=& 0; \\\\
 2 \operatorname{trace}_g A_{\nabla_{\cdot}^{\perp} H}(\cdot) + \left(p-2 + \frac{m}{2} \right) \operatorname{grad}^M |H|^2&=& 0,
\end{array}
\right.
\end{equation}
where $\Delta^{\perp}$ is the Laplacian in the normal bundle of $(M,g)$.
\end{theorem}

\begin{proof}
First, the $p$-tension field of $\mathbf{i}$ is given by
\begin{eqnarray*}
  \tau_p(\mathbf{i}) &=& |d\mathbf{i}|^{p-2} \tau(\mathbf{i}) +(p-2) |d\mathbf{i}|^{p-3} d\mathbf{i} ( \operatorname{grad}^M |d\mathbf{i}| ),
\end{eqnarray*}
since $\tau(\mathbf{i})=mH$ (see \cite{BG}, \cite{BW}), and $|d\mathbf{i}|^2=m$, we get $\tau_p(\mathbf{i})=m^{\frac{p}{2}}   H$.
Let $\{e_i,\ldots , e_m \}$ be an orthonormal frame such that $\nabla_{e_i}^Me_j =0$ at $ x \in M $ for all $i,j=1,\ldots,m $, then calculating at $x$
\begin{eqnarray*}
  \operatorname{trace}_gR^N(\tau_p(\mathbf{i}) , d\mathbf{i} ) d\mathbf{i} &=& \sum_{i=1}^{m} R^{N}(\tau_{p}(\mathbf{i}),d\mathbf{i}(e_i))d\mathbf{i}(e_i)\\
                                                                         &=& m^{\frac{p}{2}}  \sum_{i=1}^{m} R^{N}(H,e_i)e_i.
\end{eqnarray*}
By the following equation $R^N (X,Y)Z=  c \left(\langle Y,Z \rangle X - \langle X,Z \rangle Y   \right)$, with
$\langle H, e_i \rangle =0  $, for all $X,Y,Z\in\Gamma(TN(c))$ and $i=1,\ldots,m $, the last equation becomes
\begin{equation}\label{eq1}
\operatorname{trace}_{g}R^{N}(\tau_{p}(\mathbf{i}),d\mathbf{i})d\mathbf{i} = m^{\frac{p+2}{2}} c H.
\end{equation}
We compute the term $\operatorname{trace}_g(\nabla^{\mathbf{i} })^{2}\tau_{p}(\mathbf{i})$ at $x$
\begin{eqnarray}\label{eq2.1}
\sum_{i=1}^{m}\nabla^{\mathbf{i}}_{e_{i}}\nabla^{\mathbf{i}}_{e_{i}}H
  &=& \nonumber\sum_{i=1}^{m}\nabla^{\mathbf{i}}_{e_{i}}\big(-A_{H}(e_{i})+(\nabla^{\mathbf{i}}_{e_{i}}H )^{\perp}\big) \\
  &=& \nonumber -\sum_{i=1}^{m}\nabla^{M}_{e_{i}}A_{H}(e_{i})
                -\sum_{i=1}^{m}B(e_{i},A_{H}(e_{i}))\\
  & &           -\sum_{i=1}^{m}A_{(\nabla^{\mathbf{i}}_{e_{i}}H )^{\perp}}(e_{i})
                +\sum_{i=1}^{m}\big(\nabla^{\mathbf{i}}_{e_{i}}(\nabla^{\mathbf{i}}_{e_{i}}H)^{\perp}\big)^{\perp},
\end{eqnarray}
since $\langle A_{H}(X),Y \rangle=\langle B(X,Y),H \rangle$ for all $X,Y\in\Gamma(TM)$, we get
\begin{eqnarray*}
\sum_{i=1}^{m}\nabla^{M}_{e_{i}}A_{H}(e_{i})
   &=&\nonumber \sum_{i,j=1}^{m}\big\langle \nabla^{M}_{e_{i}}A_{H}(e_{i}),e_{j}\big \rangle\,e_{j} \\
   &=&\nonumber \sum_{i,j=1}^{m}e_{i}\big(\big\langle A_{H}(e_{i}),e_{j}\big \rangle\big)\,e_{j} \\
   &=&\nonumber \sum_{i,j=1}^{m} e_{i}\big(\big\langle B(e_{i},e_{j}),H\big \rangle\big)\,e_{j} \\
   &=&\nonumber \sum_{i,j=1}^{m} e_{i}\big(\big\langle \nabla^{N}_{e_{j}}e_{i},H\big \rangle\big)\,e_{j},
\end{eqnarray*}
since $\nabla^{N}_{X}\nabla^{N}_{Y}Z=
R^{N}(X,Y)Z+\nabla^{N}_{Y}\nabla^{N}_{X}Z+\nabla^{N}_{[X,Y]}Z$, for all $X,Y,Z\in\Gamma(TN(c))$, we conclude
\begin{eqnarray*}
\sum_{i=1}^{m}\nabla^{M}_{e_{i}}A_{H}(e_{i})
   &=&\nonumber \sum_{i,j=1}^{m} \big\langle \nabla^{N}_{e_{i}}\nabla^{N}_{e_{j}}e_{i},H\big \rangle\,e_{j}
                +\sum_{i,j=1}^{m} \big\langle \nabla^{N}_{e_{j}}e_{i},\nabla^{\mathbf{i}}_{e_{i}}H\big \rangle\,e_{j}\\
   &=&\nonumber \sum_{i,j=1}^{m} \big\langle R^{N}(e_{i},e_{j})e_{i},H\big \rangle\,e_{j}
                +\sum_{i,j=1}^{m} \big\langle \nabla^{N}_{e_{j}}\nabla^{N}_{e_{i}}e_{i},H\big \rangle\,e_{j}\\
   & &          +\sum_{i,j=1}^{m} \big\langle B(e_{i},e_{j}),(\nabla^{\mathbf{i}}_{e_{i}}H)^{\perp}\big \rangle\,e_{j},
\end{eqnarray*}
since $R^N (X,Y)Z=  c \left(\langle Y,Z \rangle X - \langle X,Z \rangle Y   \right)$, for all $X,Y,Z\in\Gamma(TN(c))$, we have
\begin{eqnarray}\label{eq2.2}
\sum_{i=1}^{m}\nabla^{M}_{e_{i}}A_{H}(e_{i})
&=&\nonumber \sum_{i,j=1}^{m} e_{j}\big(\big\langle \nabla^{N}_{e_{i}}e_{i},H\big \rangle\big)\,e_{j}
             -\sum_{i,j=1}^{m} \big\langle \nabla^{N}_{e_{i}}e_{i},\nabla^{\mathbf{i}}_{e_{j}}H\big \rangle\,e_{j}\\
& &\nonumber +\sum_{i,j=1}^{m}\big\langle A_{(\nabla^{\mathbf{i}}_{e_{i}}H)^{\perp}}(e_{i}),e_{j}\big \rangle\,e_{j}\\
&=&\nonumber m\sum_{j=1}^{m}e_{j}\big(\big\langle H,H\big \rangle\big)\,e_{j}
             -m\sum_{j=1}^{m}\big\langle H,\nabla^{\mathbf{i}}_{e_{j}}H\big \rangle\,e_{j}\\
& &\nonumber +\sum_{i=1}^{m}A_{(\nabla^{\mathbf{i}}_{e_{i}}H)^{\perp}}(e_{i})\\
&=&          \frac{m}{2}\sum_{j=1}^{m}e_{j}\big(\big\langle H,H\big \rangle\big)\,e_{j}
             +\sum_{i=1}^{m}A_{(\nabla^{\mathbf{i}}_{e_{i}}H)^{\perp}}(e_{i}).
\end{eqnarray}
From equations (\ref{eq2.1}) and (\ref{eq2.2}), we obtain
\begin{eqnarray}\label{eq2.3}
  \operatorname{trace}_g(\nabla^{\mathbf{i}})^{2}\tau_{p}(\mathbf{i})
  &=&\nonumber -\frac{m^{\frac{p+2}{2}} }{2}\,\operatorname{grad}^{M}|H|^{2}
     -2m^{\frac{p}{2}} \operatorname{trace}_gA_{(\nabla^{\perp}_{\cdot}H )}(\cdot)\\
  & & -m^{\frac{p}{2}} \operatorname{trace}_gB(\cdot,A_{H}(\cdot))
      +m^{\frac{p}{2}} \Delta^{\perp}H.
\end{eqnarray}
Now, we compute the term $\operatorname{trace}_g \nabla \langle  \nabla^\mathbf{i} \tau_p(\mathbf{i})  , d\mathbf{i}  \rangle d\mathbf{i}$ at $x$
\begin{eqnarray*}
\sum_{i,j=1}^m \nabla_{e_i}^\mathbf{i} \langle  \nabla_{e_j}^\mathbf{i} \, \tau_p(\mathbf{i}) , d\mathbf{i}(e_j)  \rangle d\mathbf{i}(e_i)
 &=& m^{\frac{p}{2}} \sum_{i,j=1}^m   \nabla_{e_i}^{\mathbf{i}} \langle  \nabla_{e_j}^{\mathbf{i}}   H  , e_j  \rangle e_i,
\end{eqnarray*}
by the compatibility of pull-back connection $\nabla^\mathbf{i}$ with the Riemannian metric of $N(c)$, and the definition of the mean curvature vector
field $H$ of $(M,g)$, we have
\begin{eqnarray*}
\sum_{j=1}^m \langle \nabla_{e_j}^{\mathbf{i}}H \, , e_j\rangle
&=&\sum_{j=1}^m  \big\{e_j \langle H , e_j \rangle  \, - \, \langle H , \nabla_{e_j}^{\mathbf{i}} e_j\rangle \big\}\\
&=& -\sum_{j=1}^m \langle H , B ( e_j , e_j ) \rangle  \\
&=& -m |H|^2,
\end{eqnarray*}
by the last two equations, we have the following
\begin{equation}\label{eq2}
\operatorname{trace}_g \nabla \langle  \nabla^\mathbf{i}  \tau_p(\mathbf{i})  , d\mathbf{i}  \rangle d\mathbf{i}
= - m^{\frac{p+2}{2}} \operatorname{grad}^M  |H|^2   - m^{\frac{p+4}{2}} |H|^2  H.
\end{equation}
The Theorem 2.1 followed by (\ref{eq1.7}),  (\ref{eq1}), (\ref{eq2.3}), and (\ref{eq2}).
\end{proof}
%%%%%%%%%%%%%%%%%%%%%%%%%%%%%%%%%%%%%%%%%%%%%%%%%%%%%%%%%%%%%%%%%%%%%%%%%%%%%%%%%%
%%%%%%%%%%%%%%%%%%%%%%%%%%%% Style of Definition and Remark %%%%%%%%%%%%%%%%%%%%%%
%%%%%%%%%%%%%%%%%%%%%%%%%%%% Style of Definition and Remark %%%%%%%%%%%%%%%%%%%%%%
%%%%%%%%%%%%%%%%%%%%%%%%%%%% Style of Definition and Remark %%%%%%%%%%%%%%%%%%%%%%
%%%%%%%%%%%%%%%%%%%%%%%%%%%%%%%%%%%%%%%%%%%%%%%%%%%%%%%%%%%%%%%%%%%%%%%%%%%%%%%%%%

If $p=2$ and $N=\mathbb{S}^{n}$, we arrive at the following Corollary.
\begin{corollary}
Let $M$ be a submanifold of sphere $\mathbb{S}^{n}$ of dimension $m$, then
the canonical inclusion $\mathbf{i}:M\hookrightarrow \mathbb{S}^{n}$ is biharmonic if and only if
$$\left\{
  \begin{array}{ll}
    \displaystyle\frac{m}{2}\operatorname{grad}^{M}|H|^{2}
       +2\operatorname{trace}_gA_{(\nabla^{\perp}_{\cdot}H )}(\cdot)=0, & \\\\
     \displaystyle -m\,H
      +\operatorname{trace}_gB(\cdot,A_{H}(\cdot))
      -\Delta^{\perp}H=0. &
  \end{array}
\right.$$
\end{corollary}
This result was deduced by B-Y. Chen and C. Oniciuc \cite{c}, \cite{o}.
\begin{theorem}
If $M$ is a hypersurface with nowhere zero mean curvature of $N^{m+1}(c)$, then $M$ is $p$-biharmonic if only if
  \begin{equation} \label{tan} \left\{
  \begin{array}{lll}
  -\Delta^{\perp}  H +  \big(|A|^2+m(p-2) |H|^2-mc\big)  H &=&0; \\\\
  2 A( \operatorname{grad}^M |H|) + \big(2(p-2)+m\big) |H| \operatorname{grad}^M |H| &=& 0.
  \end{array}
\right.
\end{equation}
\end{theorem}

\begin{proof}

Consider $\{e_1,\dotsc,e_m\}$ to be a local orthonormal frame field on $(M,g)$, and let $\eta$ the unit normal vector field at $(M,g)$ in $N^{m+1}(c)$.
We have
\begin{eqnarray*}
H &=& \langle H,\eta \rangle \eta \\
  &=& \frac{1}{m} \sum_{i=1}^m \langle B(e_i , e_i) , \eta  \rangle \eta \\
  &=& \frac{1}{m} \sum_{i=1}^m g(A(e_i) , e_i ) \eta  \\
  &=& \frac{1}{m} (\operatorname{trace}_g A ) \eta.
\end{eqnarray*}
Let $i=1,\ldots,m $, we compute
\begin{eqnarray*}
A_H (e_i) &=& \sum_{j=1}^m g(A_H (e_i) , e_j ) e_j \\
        &=& -\sum_{j=1}^m \langle \nabla_{e_i}^N H , e_j \rangle e_j \\
        &=& -\sum_{j=1}^m e_i\langle H , e_j \rangle e_j \, + \sum_{j=1}^m \, \langle H, B(e_i , e_j) \rangle e_j \\
        &=& \langle H,\eta \rangle\sum_{j=1}^m  \langle \eta , B(e_i , e_j)  \rangle e_j,
\end{eqnarray*}
by the last equation and the formula $\langle \eta , B(e_i , e_j)  \rangle=g(Ae_i,e_j)$, we obtain the following equation $A_H (e_i)=\langle H , \eta  \rangle A(e_i)$. So that
\begin{eqnarray}\label{eq2.4}
\sum_{i=1}^m B(e_i , A_H (e_i) )&=&\nonumber \sum_{i=1}^m B(e_i , \langle H , \eta \rangle A (e_i) ) \\
                 &=&\nonumber \langle H,\eta \rangle \sum_{i=1}^m B ( e_i , A (e_i)) \\
                 &=&\nonumber \langle H , \eta  \rangle \sum_{i=1}^m g( A(e_i) , A(e_i))\eta \\
                 &=& |A|^2  H.
\end{eqnarray}
In the same way, with $\eta=H/|H|$, we find that
\begin{eqnarray}\label{eq2.5}
\sum_{i=1}^mA_{\nabla_{e_i}^\perp H} (e_i)
                              &=&\nonumber \sum_{i,j=1}^m\langle A_{\nabla_{e_i}^\perp H} (e_i) , e_j  \rangle e_j \\
                              &=&\nonumber -\sum_{i,j=1}^m\langle \nabla_{e_i}^N \nabla_{e_i}^\perp H , e_j \rangle e_j \\
                              &=&\nonumber - \sum_{i,j=1}^m\langle e_i\langle H, \eta \rangle\nabla_{e_i}^N \eta , e_j \rangle e_j\\
                              &=& A(\operatorname{grad}^M |H|).
\end{eqnarray}
The Theorem 2.3 followed by equations (\ref{eq2.4}),  (\ref{eq2.5}), and Theorem 2.1.
\end{proof}

\begin{corollary}
$(i)$ A submanifold $M$ with parallel mean curvature vector field in $N^n(c)$ is $p$-biharmonic if and only if
\begin{equation}
\operatorname{trace}_g B(\cdot , A_H (\cdot) ) = m\big(c-(p-2) |H|^2\big)  H ,
\end{equation}
$(ii)$ A hypersurface $M$ of constant non-zero mean curvature in $N^{m+1}(c)$ is proper  $p$-biharmonic if and only if
\begin{equation}
 |A|^2  = mc- m(p-2) |H|^2.
\end{equation}
\end{corollary}
\begin{example}
We consider the hypersurface
$$ \mathbb{S}^m(a)  = \big\{  (x^1, \cdots , x^m , x^{m+1} , b ) \, \in \mathbb{R}^{m+2} \, : \, \sum_{i=1}^{m+1} (x^i)^2 = a^2 \big\} \subset  \mathbb{S}^{m+1}, $$
where $a^2+b^2=1$. We have
$$\eta = \frac{1}{r} (x^1, \cdots , x^{m+1} , -\frac{a^2}{b}) , $$
with $r^2 = \frac{a^2}{b^2}$ $(r> 0) $, is a unit section in the normal bundle of $\mathbb{S}^m(a)$ in $\mathbb{S}^{m+1}$.\\ Let $X \in  \Gamma ( T \mathbb{S}^m (a))$,
we compute
$$\nabla_X^{\mathbb{S}^{m+1}} \eta = \frac{1}{r} \nabla_X^{\mathbb{R}^{m+2}} ( x^1 ,\cdots , x^{m+1} , -\frac{a^2}{b}) = \frac{1}{r} X.$$
Thus, $ \nabla^{\perp} \eta = 0 $ and $A = - \frac{1}{r} Id$.
This implies that $H = - \frac{1}{r} \eta $, and so $\mathbb{S}^m(a) $ has constant mean curvature $|H| = \frac{1}{r} $ in $\mathbb{S}^{m+1}$.
Since $|A|^2 = \frac{m}{r^2}$, according to Corollary 2.4. we conclude that $\mathbb{S}^m( a) $ is proper $p$-biharmonic in $\mathbb{S}^{m+1}$ if and only if $p=1/b^2$.
\end{example}

\section{Stress $p$-bienergy tensors}
Let $\varphi:(M,g)\rightarrow (N,h)$ be a smooth map between two Riemannian manifolds and  $p\geq2$. Consider a smooth one-parameter variation of the metric $g$, i.e. a smooth family of metrics $(g_{t})$ $(-\epsilon<t<\epsilon)$ such that $g_{0}=g$, write $\delta=\frac{\partial}{\partial t}\:\Big|_{t=0}$, then
$\delta g\in\Gamma(\odot^{2}T^{*}M)$ is a symmetric $2$-covariant tensor field on $M$ (see \cite{BW}).
Take local coordinates $(x^{i})$ on $M$, and write the metric on $M$ in the usual way as
$g_{t}=g_{ij}(t,x)\,dx^{i}\, dx^{j}$, we now compute
\begin{equation}\label{eq3.2}
\frac{d}{dt}E_{2,p}(\varphi;D)\Big|_{t=0}=\frac{1}{2}\int_{D}\delta(|\tau_{p}(\varphi)|^{2})v_{g}+
\frac{1}{2}\int_{D}|\tau_{p}(\varphi)|^{2}\delta( v_{g_{t}}).
\end{equation}
The calculation of the first term breaks down in three lemmas.
\begin{lemma}\label{lemma1}
The vector field $\xi=( \operatorname{div}^{M}\delta g)^{\sharp}-\frac{1}{2}\operatorname{grad}^{M}(\operatorname{trace}\,\delta g)$ satisfies
\begin{eqnarray*}
   \delta(|\tau_{p}(\varphi)|^{2})
   &=&-(p-2)|d\varphi|^{p-4}\langle \varphi^{*}h,\delta g \rangle h(\tau(\varphi),\tau_{p}(\varphi))\\
   &&-2|d\varphi|^{p-2}\langle  h(\nabla d\varphi,\tau_{p}(\varphi)),\delta g \rangle
       -2|d\varphi|^{p-2}h(d\varphi(\xi),\tau_{p}(\varphi))\\
   &&-(p-2)(p-4)|d\varphi|^{p-5}\langle \varphi^{*}h,\delta g\rangle h(d\varphi(\operatorname{grad}^{M}|d\varphi|),\tau_{p}(\varphi))\\
   &&-2(p-2)|d\varphi|^{p-3}\langle d|d\varphi|\odot h(d\varphi,\tau_{p}(\varphi)),\delta g \rangle\\
   &&-(p-2)|d\varphi|^{p-4}h(d\varphi(\operatorname{grad}^{M}\langle \varphi^{*}h,\delta g\rangle),\tau_{p}(\varphi)),
\end{eqnarray*}
where $\varphi^{*}h$ is the pull-back of the metric $h$, and $\langle\,,\,\rangle$ is the induced Riemannian metric on $\otimes^{2}T^{*}M$.
\end{lemma}
\begin{proof}
In local coordinates $(x^{i})$ on $M$ and $(y^{\alpha})$ on $N$, we have
\begin{equation}\label{eq3.3}
\delta(|\tau_{p}(\varphi)|^{2})=
\delta\big(\tau_{p}(\varphi)^{\alpha}\tau_{p}(\varphi)^{\beta}h_{\alpha\beta}\big)=
2\delta(\tau_{p}(\varphi)^{\alpha})\tau_{p}(\varphi)^{\beta}h_{\alpha\beta}.
\end{equation}
By the definition of $\tau_{p}(\varphi)$ we get
\begin{eqnarray}\label{eq3.4}
\delta(\tau_{p}(\varphi)^{\alpha})
   &=& \nonumber\delta\big(|d\varphi|^{p-2}\tau(\varphi)^{\alpha}+\theta^{\alpha}\big) \\
   &=& \delta(|d\varphi|^{p-2})\tau(\varphi)^{\alpha}+|d\varphi|^{p-2}\delta(\tau(\varphi)^{\alpha})+\delta(\theta^{\alpha}).
\end{eqnarray}
where
$
\tau(\varphi)^{\alpha}=g^{ij}\big(\varphi^{\alpha}_{i,j}
+^{N}\Gamma_{\mu\sigma}^{\alpha}\varphi_{i}^{\mu}\varphi_{j}^{\sigma}
-^{M}\Gamma_{ij}^{k}\varphi^{\alpha}_{k}\big)
$ is the component of the tension field $\tau(\varphi)$, and
$\theta^{\alpha}=(p-2)|d\varphi|^{p-3}g^{ij}|d\varphi|_i\varphi_{j}^{\alpha}.$\\
The first term in the right-hand side of (\ref{eq3.4}) is given by
\begin{eqnarray}\label{eq3.5}
\delta(|d\varphi|^{p-2})\,\tau(\varphi)^{\alpha}
&=&\nonumber(p-2)|d\varphi|^{p-4}\delta(\frac{|d\varphi|^{2}}{2})\tau(\varphi)^{\alpha}\\
&=&-\frac{p-2}{2}|d\varphi|^{p-4}\langle \varphi^{*}h,\delta g \rangle\tau(\varphi)^{\alpha}.
\end{eqnarray}
The second term on the right-hand side of (\ref{eq3.4}) is (see \cite{LMO})
\begin{equation}\label{eq3.6}
|d\varphi|^{p-2}\delta(\tau(\varphi)^{\alpha})=
-|d\varphi|^{p-2}g^{ai}g^{bj}\delta(g_{ab}) (\nabla d\varphi)_{ij}^{\alpha}-|d\varphi|^{p-2}\xi^{k}\varphi_{k}^{\alpha},
\end{equation}
Now, we compute the third term on the right-hand side of (\ref{eq3.4})
\begin{eqnarray}\label{eq3.7}
\delta(\theta^{\alpha})
   &=&\nonumber (p-2)(p-3)|d\varphi|^{p-5}\delta(\frac{|d\varphi|^{2}}{2})g^{ij}|d\varphi|_i\varphi_{j}^{\alpha}\\
   &&\nonumber+(p-2)|d\varphi|^{p-3}\delta(g^{ij})|d\varphi|_i\varphi_{j}^{\alpha}\\
   &&+(p-2)|d\varphi|^{p-3}g^{ij}\delta(|d\varphi|_i)\varphi_{j}^{\alpha}.
\end{eqnarray}
By using $\delta(\frac{|d\varphi|^{2}}{2})=-\frac{1}{2}\langle \varphi^{*}h,\delta g\rangle$ with
$\delta(|d\varphi|_i)=(\delta(|d\varphi|))_i$, the equation (\ref{eq3.7}) becomes
\begin{eqnarray}\label{eq3.8}
\delta(\theta^{\alpha})
   &=&\nonumber -\frac{(p-2)(p-3)}{2}|d\varphi|^{p-5}\langle \varphi^{*}h,\delta g\rangle g^{ij}|d\varphi|_i\varphi_{j}^{\alpha}\\
   &&\nonumber+(p-2)|d\varphi|^{p-3}\delta(g^{ij})|d\varphi|_i\varphi_{j}^{\alpha}\\
   &&\nonumber-\frac{p-2}{2}|d\varphi|^{p-4}g^{ij}\langle \varphi^{*}h,\delta g\rangle_i\varphi_{j}^{\alpha}\\
   &&+\frac{p-2}{2}|d\varphi|^{p-5}g^{ij}|d\varphi|_i\langle \varphi^{*}h,\delta g\rangle\varphi_{j}^{\alpha}.
\end{eqnarray}
Note that
\begin{eqnarray}\label{eq3.9}
2\delta(|d\varphi|^{p-2})\tau(\varphi)^{\alpha}\tau_{p}(\varphi)^{\beta}h_{\alpha\beta}
&=&\nonumber-(p-2)|d\varphi|^{p-4}\langle \varphi^{*}h,\delta g \rangle\tau(\varphi)^{\alpha}\tau_{p}(\varphi)^{\beta}h_{\alpha\beta}\\
&=&\nonumber-(p-2)|d\varphi|^{p-4}\langle \varphi^{*}h,\delta g \rangle h(\tau(\varphi),\tau_{p}(\varphi)),\\
\end{eqnarray}
\begin{eqnarray}\label{eq3.10}
2
|d\varphi|^{p-2}\delta(\tau(\varphi)^{\alpha})
\tau_{p}(\varphi)^{\beta}h_{\alpha\beta}
   &=& \nonumber -2|d\varphi|^{p-2}g^{ai}g^{bj}\delta(g_{ab}) (\nabla d\varphi)_{ij}^{\alpha}\tau_{p}(\varphi)^{\beta}h_{\alpha\beta}\\
   & & \nonumber -2|d\varphi|^{p-2}\xi^{k}\varphi_{k}^{\alpha}\tau_{p}(\varphi)^{\beta}h_{\alpha\beta}\\
   &=& \nonumber-2|d\varphi|^{p-2}\langle  h(\nabla d\varphi,\tau_{p}(\varphi)),\delta g \rangle\\
   & &    -2|d\varphi|^{p-2}h(d\varphi(\xi),\tau_{p}(\varphi)),
\end{eqnarray}
and the following
\begin{eqnarray}\label{eq3.11}
2
\delta(\theta^{\alpha})
\tau_{p}(\varphi)^{\beta}h_{\alpha\beta}
 &=&\nonumber -(p-2)(p-3)|d\varphi|^{p-5}\langle \varphi^{*}h,\delta g\rangle g^{ij}|d\varphi|_i\varphi_{j}^{\alpha}\tau_{p}(\varphi)^{\beta}h_{\alpha\beta}\\
   &&\nonumber+2(p-2)|d\varphi|^{p-3}\delta(g^{ij})|d\varphi|_i\varphi_{j}^{\alpha}\tau_{p}(\varphi)^{\beta}h_{\alpha\beta}\\
   &&\nonumber-(p-2)|d\varphi|^{p-4}g^{ij}\langle \varphi^{*}h,\delta g\rangle_i\varphi_{j}^{\alpha}\tau_{p}(\varphi)^{\beta}h_{\alpha\beta}\\
   &&\nonumber+(p-2)|d\varphi|^{p-5}g^{ij}|d\varphi|_i\langle \varphi^{*}h,\delta g\rangle\varphi_{j}^{\alpha}\tau_{p}(\varphi)^{\beta}h_{\alpha\beta}\\
   &=&\nonumber -(p-2)(p-3)|d\varphi|^{p-5}\langle \varphi^{*}h,\delta g\rangle h(d\varphi(\operatorname{grad}^{M}|d\varphi|),\tau_{p}(\varphi))\\
   &&\nonumber-2(p-2)|d\varphi|^{p-3}\langle d|d\varphi|\odot h(d\varphi,\tau_{p}(\varphi)),\delta g \rangle\\
   &&\nonumber-(p-2)|d\varphi|^{p-4}h(d\varphi(\operatorname{grad}^{M}\langle \varphi^{*}h,\delta g\rangle),\tau_{p}(\varphi))\\
   &&+(p-2)|d\varphi|^{p-5}\langle \varphi^{*}h,\delta g\rangle h(d\varphi(\operatorname{grad}^{M}|d\varphi|)  ,\tau_{p}(\varphi)).
\end{eqnarray}
Substituting (\ref{eq3.4}), (\ref{eq3.9}), (\ref{eq3.10}) and (\ref{eq3.11}) in (\ref{eq3.3}), the Lemma \ref{lemma1} follows.
\end{proof}
\begin{lemma}[\cite{DC}] \label{lemma2} Let $D$ be a compact domain of $M$. Then
\begin{eqnarray*}
  \int_{D}|d\varphi|^{p-2}h(d\varphi(\xi),\tau_{p}(\varphi))v_{g} &=&
  \int_{D}\big\langle -\operatorname{sym}\big(\nabla |d\varphi|^{p-2}h(d\varphi,\tau_{p}(\varphi))\big) \\
   &&+\frac{1}{2}\operatorname{div}^{M}\big(|d\varphi|^{p-2}h(d\varphi,\tau_{p}(\varphi))
^{\sharp}\big)g,\delta g\big\rangle \,v_{g}.
\end{eqnarray*}
\end{lemma}
\begin{lemma}\label{lemma3} We set $\omega=|d\varphi|^{p-4}h(d\varphi,\tau_{p}(\varphi))$. Then
\begin{eqnarray*}
  -\int_{D}|d\varphi|^{p-4}h(d\varphi(\operatorname{grad}^{M}\langle \varphi^{*}h,\delta g \rangle),\tau_{p}(\varphi))\,v_{g} &=&
  \int_{D}\langle \varphi^{*}h,\delta g \rangle\operatorname{div}\omega \,v_{g}.
\end{eqnarray*}
\end{lemma}
\begin{proof}
Note that
$$\operatorname{div}(\langle \varphi^{*}h,\delta g \rangle\omega)
=\langle \varphi^{*}h,\delta g \rangle\operatorname{div}\omega
+\omega(\operatorname{grad}^{M}\langle \varphi^{*}h,\delta g \rangle),$$
and consider the divergence Theorem, Lemma \ref{lemma3} follows.
\end{proof}
\begin{theorem}\label{Theorem1}
Let $\varphi:(M,g)\rightarrow(N,h)$ be a smooth map such that $|d\varphi|_x\neq0$ for all $x\in M$,
and let $\{g_{t}\}$ a one parameter variation of $g$. Then
$$\frac{d}{dt}E_{2,p}(\varphi;D)\Big|_{t=0}=\frac{1}{2}\int_{D}\langle S_{2,p}(\varphi),\delta g \rangle\,v_{g},$$
where $S_{2,p}(\varphi)\in\Gamma(\odot^{2}T^{*}M)$ is given by
\begin{eqnarray*}
  S_{2,p}(\varphi) (X,Y)
   &=& -\frac{1}{2}|\tau_{p}(\varphi)|^{2}g(X,Y)
       -|d\varphi|^{p-2}\langle d\varphi,\nabla^{\varphi} \tau_{p}(\varphi) \rangle g(X,Y)\\
   & & +|d\varphi|^{p-2}h(d\varphi(X),\nabla_{Y}^{\varphi}\tau_{p}(\varphi))
       +|d\varphi|^{p-2}h(d\varphi(Y),\nabla_{X}^{\varphi}\tau_{p}(\varphi))\\
   & & +(p-2)|d\varphi|^{p-4}\langle d\varphi,\nabla^{\varphi} \tau_{p}(\varphi) \rangle h(d\varphi(X),d\varphi(Y)).
\end{eqnarray*}
$S_ {2,p}(\varphi) $ is called the stress $p$-bienergy tensor of $\varphi$.
\end{theorem}

\begin{proof}
By using $\delta(v_{g_{t}})=\frac{1}{2}\langle g,\delta g\rangle v_{g}$ (see \cite{BW}),
 Lemmas \ref{lemma1}, \ref{lemma2}, and \ref{lemma3}, the equation (\ref{eq3.2}) becomes
\begin{eqnarray}\label{eq3.13}
   S_ {2,f}(\varphi)
   &=&\nonumber-(p-2)|d\varphi|^{p-4}h(\tau(\varphi),\tau_{p}(\varphi))\varphi^{*}h\\
   &&\nonumber-2|d\varphi|^{p-2} h(\nabla d\varphi,\tau_{p}(\varphi))
       +2\operatorname{sym}\big(\nabla |d\varphi|^{p-2}h(d\varphi,\tau_{p}(\varphi))\big) \\
   &&\nonumber-\operatorname{div}^{M}\big(|d\varphi|^{p-2}h(d\varphi,\tau_{p}(\varphi))
^{\sharp}\big)g   \\
   &&\nonumber-(p-2)(p-4)|d\varphi|^{p-5} h(d\varphi(\operatorname{grad}^{M}|d\varphi|),\tau_{p}(\varphi)) \varphi^{*}h\\
   &&\nonumber-2(p-2)|d\varphi|^{p-3}d|d\varphi|\odot h(d\varphi,\tau_{p}(\varphi))\\
   &&+(p-2)\operatorname{div}^{M}\big[|d\varphi|^{p-4}h(d\varphi,\tau_{p}(\varphi))\big]\varphi^{*}h
   +\frac{1}{2}|\tau_{p}(\varphi)|^2 g.
\end{eqnarray}
Note that, for all $X,Y\in\Gamma(TM)$, we have
\begin{eqnarray}\label{eq3.14}
 2\operatorname{sym}\big(\nabla |d\varphi|^{p-2}h(d\varphi,\tau_{p}(\varphi))\big)(X,Y)
 &=&\nonumber 2|d\varphi|^{p-2}h(\nabla d\varphi(X,Y),\tau_{p}(\varphi))  \\
 & &\nonumber +|d\varphi|^{p-2}h(d\varphi(X),\nabla^{\varphi}_{Y}\tau_{p}(\varphi))\\
 & &\nonumber +|d\varphi|^{p-2}h(d\varphi(Y),\nabla^{\varphi}_{X}\tau_{p}(\varphi))\\
 & &\nonumber +X(|d\varphi|^{p-2})h(d\varphi(Y),\tau_{p}(\varphi))   \\
 & &\nonumber +Y(|d\varphi|^{p-2})h(d\varphi(X),\tau_{p}(\varphi)),\\
\end{eqnarray}
and the following formula
\begin{eqnarray}\label{eq3.15}
 - 2 d|d\varphi|\odot h(d\varphi,\tau_{p}(\varphi))(X,Y)
   &=&\nonumber -X(|d\varphi|)h(d\varphi(Y),\tau_{p}(\varphi))  \\
   & & -Y(|d\varphi|)h(d\varphi(X),\tau_{p}(\varphi)).
\end{eqnarray}
Calculating in a normal frame at $x$, we have
\begin{eqnarray}\label{eq3.16}
  \operatorname{div}^{M}\big(|d\varphi|^{p-2}h(d\varphi,\tau_{p}(\varphi))^{\sharp}\big)
  &=&\nonumber \sum_{i=1}^me_{i}(g( |d\varphi|^{p-2}h(d\varphi,\tau_{p}(\varphi))^{\sharp},e_{i}))\\
  &=&\nonumber \sum_{i=1}^me_{i}(|d\varphi|^{p-2}h(d\varphi(e_{i}),\tau_{p}(\varphi)))\\
  &=&\nonumber \sum_{i=1}^me_{i}( |d\varphi|^{p-2})h(d\varphi(e_{i}),\tau_{p}(\varphi))\\
  && \nonumber   +\sum_{i=1}^m|d\varphi|^{p-2}h(\nabla^{\varphi}_{e_{i}}d\varphi(e_{i}),\tau_{p}(\varphi))\\
  & &\nonumber +\sum_{i=1}^m|d\varphi|^{p-2}h(d\varphi(e_{i}),\nabla^{\varphi}_{e_{i}}\tau_{p}(\varphi))\\
  &=&\nonumber (p-2)|d\varphi|^{p-3}h(d\varphi(\operatorname{grad}^{M}|d\varphi|),\tau_{p}(\varphi))\\
    &&\nonumber  +|d\varphi|^{p-2}h(\tau(\varphi),\tau_{p}(\varphi))\\
  & & +|d\varphi|^{p-2}\langle d\varphi,\nabla^{\varphi}\tau_{p}(\varphi) \rangle.
\end{eqnarray}
From the definition of $\tau_{p}(\varphi)$, and equation (\ref{eq3.16}), we get
\begin{eqnarray}\label{eq3.17}
\operatorname{div}^{M}\big(|d\varphi|^{p-2}h(d\varphi,\tau_{p}(\varphi))^{\sharp}\big)
=|\tau_{p}(\varphi)|^2
+|d\varphi|^{p-2}\langle d\varphi,\nabla^{\varphi}\tau_{p}(\varphi) \rangle.
\end{eqnarray}
With the same method of (\ref{eq3.16}), we find that
\begin{eqnarray}\label{eq3.18}
  \operatorname{div}^{M}\big(|d\varphi|^{p-4}h(d\varphi,\tau_{p}(\varphi))\big)
  &=&\nonumber (p-4)|d\varphi|^{p-5}h(d\varphi(\operatorname{grad}^{M}|d\varphi|),\tau_{p}(\varphi))\\
    && \nonumber +|d\varphi|^{p-4}h(\tau(\varphi),\tau_{p}(\varphi))\\
  & & +|d\varphi|^{p-4}\langle d\varphi,\nabla^{\varphi}\tau_{p}(\varphi) \rangle.
\end{eqnarray}
Substituting (\ref{eq3.14}), (\ref{eq3.15}), (\ref{eq3.17}) and (\ref{eq3.18}) in (\ref{eq3.13}), the Theorem \ref{Theorem1} follows.
\end{proof}
By using the definition of divergence for symmetric $(0,2)$-tensors (see \cite{BW}, \cite{DC})
we have the following result.
\begin{theorem}
Let $\varphi:(M,g)\rightarrow(N,h)$ be a smooth map such that $|d\varphi|_x\neq0$ for all $x\in M$. Then
\begin{eqnarray*}
  \operatorname{div}^{M} S_{2,p}(\varphi) (X)&=&-h(\tau_{2,p}(\varphi),d\varphi(X)),\quad\forall X\in\Gamma(TM).
\end{eqnarray*}
\end{theorem}

\begin{remark}
When $p=2$, we have $S_{2,p}(\varphi)=S_{2}(\varphi)$, where $S_{2}(\varphi)$ is stress bienergy tensor
in \cite{LMO}.
\end{remark}

\begin{corollary}
Let $\varphi:(M,g)\rightarrow(N,h)$ be a smooth map. (1) Then $S_{2,m}(\varphi)=0$ implies that $\varphi$
is $m$-harmonic, where $m=\dim M$. (2) If $M$ is compact without
boundary, and $p\neq\frac{m}{2}$. Then $S_{2,p}(\varphi)=0$ implies $\varphi$
is $p$-harmonic.
\end{corollary}
\begin{proof}
Let $\{e_{i}\}$ be an orthonormal frame on $(M,g)$. (1) We have
\begin{eqnarray*}
  0=\sum_{i=1}^mS_{2,p}(\varphi) (e_{i},e_{i})
   &=& -\frac{m}{2}|\tau_{p}(\varphi)|^{2}
       +(p-m)|d\varphi|^{p-2}\langle d\varphi,\nabla^{\varphi} \tau_{p}(\varphi) \rangle.
\end{eqnarray*}
For $p=m$, the last equation becomes $-\frac{m}{2}|\tau_{m}(\varphi)|^{2}=0$. So $\varphi$
is $m$-harmonic map.
(2) We set $\theta(X)=h(|d\varphi|^{p-2}d\varphi(X), \tau_{p}(\varphi))$, for all $X\in\Gamma(TM)$.
The trace of $S_{2,p}(\varphi)$ gives the equality
\begin{eqnarray*}
  0=\sum_{i=1}^mS_{2,p}(\varphi) (e_{i},e_{i})
   &=& (\frac{m}{2}-p)|\tau_{p}(\varphi)|^{2}
       +(p-m)\operatorname{div}^{M}\theta.
\end{eqnarray*}
By using the Green Theorem, we get
\begin{eqnarray*}
  (\frac{m}{2}-p)\int_M|\tau_{p}(\varphi)|^{2}v^g=0.
\end{eqnarray*}
Since $p\neq\frac{m}{2}$, we obtain $|\tau_{p}(\varphi)|^2=0$, that is $\varphi$
is $p$-harmonic map.
\end{proof}

%%% References

\EditInfo{June 01, 2021}{December 03, 2021}{Haizhong Li}

\end{paper}